\documentclass[a4paper,10pt]{article}
\usepackage{amssymb}
\usepackage{amsmath,amsfonts,amsthm,amssymb}
\usepackage[dvips]{graphics}
\usepackage{epsfig}
\usepackage{indentfirst}
\usepackage{color}
\allowdisplaybreaks

\newtheoremstyle{theorem}
  {10pt}          
  {10pt}  
  {\sl}  
 {}
  {\bf}  
  {. }    
  { }    
  {}     
\theoremstyle{theorem}

\newtheorem{theorem}{Theorem}[section]

 \newtheorem{remark}{Remark}[section]
  
\numberwithin{equation}{section}

\newtheoremstyle{defi}
  {10pt}          
  {10pt}  
  {\rm}  
  {}  
  {\bf}  
  {. }    
  { }    
  {}     
\theoremstyle{defi}

\begin{document}
\baselineskip = 13.5pt

\title{\bf On a singular limit for the compressible rotating Euler system  }

\author{\v{S}\'{a}rka Ne\v{c}asov\'{a}$^1$\thanks{The research of \v S.N. leading to these results has received funding from
the Czech Sciences Foundation (GA\v CR),
P201-16-032308
 and RVO 67985840. Final version of the paper was made under support the Czech Sciences Foundation (GA\v CR),  GA19-04243S. }
\footnote{Email: matus@math.cas.cz} \ \ \  Tong Tang$^{2,1}$\thanks{The research of T.T. is supported by the NSFC Grant No. 11801138.} \footnote{Email: tt0507010156@126.com}\\
{\small  1. Institute of Mathematics of the Academy of Sciences of the Czech Republic,} \\
{\small \v Zitn\' a 25, 11567, Praha 1, Czech Republic}\\
{\small 2. Department of Mathematics, College of Sciences,}\\
{\small Hohai University, Nanjing 210098, P.R. China}\\
\date{}}
\maketitle
\begin{abstract}
The work addresses a singular limit for a rotating compressible Euler system in the low Mach number and low Rossby number regime. Based on the concept of dissipative measure-valued solution, the quasi-geostrophic system is identified as the limit problem in the case of ill-prepared initial data. The ill-prepared initial data will cause rapidly oscillating acoustic waves. Using dispersive estimates of Strichartz type, the effect of the acoustic waves in the asymptotic limit is eliminated.
\vspace{0.5cm}


\medskip

\end{abstract}
\maketitle
\section{Introduction}\setcounter{equation}{0}
Earth's graceful rotation is an unignorable factor at geophysical fluids models. These models play an important role in the analysis of complex Earth phenomena in meteorology, geophysical and astrophysics. In order to describe the effect of rotation, people introduce two factors: Coriolis acceleration and centrifugal acceleration. In many real world applications, the action of centrifugal force is neglected, as it is in equilibrium with stratification caused by the gravity of the Earth. Under the above assumptions, we consider the following scaled Euler equations in an infinite slap $\Omega=\mathbb{R}^2\times(0,1)$:
\begin{eqnarray}
\left\{
\begin{array}{llll}  \partial_{t}\rho+\text{div}(\rho \mathbf{u})=0, \\
\partial_t(\rho \mathbf{u})+\textrm{div}(\rho\mathbf{u}\otimes\mathbf{u})+\frac{1}{Ma^2}\nabla_x p(\rho)+\frac{1}{Ro}\rho(\omega\times\mathbf{u})=0,\label{1.3}
\end{array}\right.
\end{eqnarray}
where the unknown fields $\rho=\rho(t,x)$ and $\mathbf{u}=\mathbf{u}(t,x)$
represent the density and the velocity of an inviscid compressible fluid, $\omega=(0,0,1)$ is the rotation axis. The Mach number Ma, proportional to the characteristic velocity field divided by the sound speed, and the Rossby number Ro, defined as the ratio of the displacement due to Coriolis forces, play the role of singular (small) parameters. The symbol $p = p(\rho)$ denotes the barotropic pressure (assumptions on the pressure see (\ref{pres})). The system is supplemented by the far field conditions
\begin{eqnarray}
\mathbf{u}\rightarrow0,\hspace{10pt}\rho\rightarrow\overline{\rho}, \hspace{5pt} \text{as}\hspace{3pt} |x|\rightarrow\infty,
\ \mbox{where} \hspace{5pt}\overline{\rho}>0,\label{1.2}
\end{eqnarray}
and boundary condition
\begin{eqnarray}
\mathbf{u}\cdot\mathbf{n}|_{\partial\Omega}=0,\label{1.3}
\end{eqnarray}
where $\mathbf{n}$ is outer normal vector to $\partial\Omega$.

From modeling of geophysical fluids, the value of Mach number and Rossby number can be considered very small. It is well known that the compressible fluid flow becomes incompressible in the low Mach number limit, as the density distribution is constant and the velocity field becomes solenoidal. On the other hand, low Rossby number corresponds to fast rotation and the fast rotating fluids will lead to the so-called Taylor-Proudman columns phenomena. Therefore, it is interesting to observe the phenomenon if the two effects take place simultaneously. In this paper, we address the problem of the double limit for $Ma = Ro = \epsilon$. Let $\rho = \rho_\epsilon$, $\mathbf{u} = \mathbf{u}_\epsilon$,
the system (1.1) takes the form
\begin{eqnarray}
\left\{
\begin{array}{llll}  \partial_{t}\rho_\epsilon+\text{div}(\rho_\epsilon \mathbf{u}_\epsilon)=0, \\
\partial_t(\rho_\epsilon\mathbf{u}_\epsilon)+\textrm{div}(\rho_\epsilon\mathbf{u}_\epsilon\otimes\mathbf{u}_\epsilon)+\frac{1}{\epsilon^2}\nabla_x p(\rho_\epsilon)+\frac{1}{\epsilon}\rho_\epsilon(\omega\times\mathbf{u}_\epsilon)=0.\label{1.4}
\end{array}\right.
\end{eqnarray}
Our goal is to study the singular limit $\epsilon \to 0$ at the case of the \emph{ill--prepared} initial data for the scaled system \eqref{1.4}. The definition of \emph{ill--prepared} initial data will be introduced in Section \ref{ILP}. Supposing we know that in the corresponding spaces,
\begin{eqnarray*}
\rho_\epsilon^{(1)}=\frac{\rho_\epsilon-\overline{\rho}}{\epsilon}\rightarrow q,\hspace{5pt} \mathbf{u}_\epsilon\rightarrow\mathbf{v},
\end{eqnarray*}
we can find that $q$ and $\mathbf{v}$ satisfy the following equations:
\begin{align}
&\omega\times\mathbf{v}+\frac{p'(\overline{\rho})}{\overline{\rho}}\nabla_xq=0,\label{1.5}\\
&\partial_t(\Delta_hq-\frac{1}{p'(\overline{\rho})}q)+\nabla_h^{\perp}q\cdot\nabla_h(\Delta_hq)=0.\label{1.6}
\end{align}
Equations \eqref{1.5}\eqref{1.6} can be interpreted as a kind of stream function, according to physicists, named as quasi-physical flows \cite{z}. For non-rotating compressible Euler fluids, a great number of well-posedness results have been obtained. However, some classical literatures show that smooth solutions of the Euler system will exhibit blow-up phenomena in a finite time no matter how smooth or small the initial data are. Therefore, it seems more appropriate to consider a suitable class of admissible weak solutions to \eqref{1.4}. By admissible we mean that solutions
will satisfy some form of the energy balance. The need for global admissible solutions of the Euler system leads to the concept of more general
\emph{dissipative measure--valued \emph{(DMV)} solutions} introduced in the context of the full Euler system in \cite{j1,j2}.

The measure-valued solutions to hyperbolic conservations laws were introduced by DiPerna \cite{Di}.
He used Young measures to pass to the artificial viscosity limit. In the case of the incompressible Euler equations, DiPerna and Majda \cite{DM}
proved the global existence of measure-valued solutions for any initial data with finite energy.
They introduced generalized Young measures to take into account oscillations and concentrations. Further,
 the existences of measure-valued solutions were shown for further models of fluids, e.g. compressible Euler and Navier-Stokes equations \cite {Ne, KW}.
The measure-valued solution to the non-Newtonian case was proved by Novotn\' y and Ne\v casov\' a \cite{NoN}.
The generalization was given by Alibert and Bouchitt\' e \cite{AB}.
The weak-strong uniqueness for generalized measure-valued solutions of isentropic Newtonian Euler equations were proved in \cite{g}.
Inspired by previous results, the concept of dissipative measure-valued solution was finally applied to the barotropic compressible Navier-Stokes system \cite {e7}.

The reader may consult \cite{e5,e6,KW,m,Ne} for applications of the theory of \emph{(DMV)} solutions in fluid mechanics or their counterparts \cite{DeStTz,D} in  other areas of mathematical physics.

Let us discuss the main differences between  weak solutions and \emph{(DMV)} solutions. First important advanatge of \emph{(DMV)} solution is that {\it  DMV solutions to the compressible Euler system exist globally in time}. Secondly {\it (DMV) solutions convergence to the limit system holds for any ill-prepared initial data}, which in both case are not valid for weak solutions.

Due to the above fascinate advantage, there are some new results concerning singular limits in the context of measure--valued solutions. The low Mach number limit was studied in \cite{e6}, where it is shown that \emph{(DMV)} solutions approach the smooth solutions of incompressible Euler system both for well-prepared and ill-prepared data. Moreover, the singular limit of compressible Euler system in the low Mach number and strong stratification regime for the ill-prepared data was identified, see \cite{b}. However, to the best of our knowledge, compared with non-rotating case, there is a few results concerning on the singular limit of rotating compressible Euler system no matter weak solutions or strong solutions. Nilasis \cite{nc} proved the singular limit of a rotating compressible Euler system with stratification at the case of well-prepared initial data. Our goal is to consider the asymptotic limit of \emph{(DMV)} solutions to the compressible Euler equations with ill-prepared initial data. We prove it converges to the strong solutions of quasi-physical flows. Moreover, we should emphasize that boundary conditions in this paper can be replaced by the periodical conditions in $x_3$ direction with a few changes. All the choices of boundary conditions prevent the flow from creating a viscous boundary layer. The periodic domain with well-prepared initial data was considered by Feireisl et al., \cite{e6}. If the whole domain is torus $T^3$, it is difficult to obtain the analysis of acoustic waves at the case of ill-prepared initial data. It seems interesting to compare the results of the present paper with those obtained in \cite{e6}. The analysis in \cite{e6} leans that the \emph{(DMV)} solutions of Euler system will converge to incompressible Euler system. Moreover, there is obvious difference about acoustic wave analysis between rotating and non-rotating case. The extension of the results of \cite{e6} to the rotating Euler system is therefore not straightforward. Last but not least, we should emphasis that there are huge results about rotating Navier-Stokes system such as \cite{ca,ch,e2,e4,no}.

The paper is organized as follows. In Section 2, we introduce the dissipative measure solutions, relative energy and the other necessary material. In Section 3, we state our main theorem. Section 4 is devoted to deriving uniform bounds of the Euler system independent of $\epsilon$. In Section 5, we perform the necessary analysis of the acoustic waves. The proof of the main theorem is completed in Section 6.

\section{Preliminaries}

First let us observe that it is more convenient to rewrite the Euler system in terms of the conservative variables $\rho$, $\mathbf{m} = \rho \mathbf{u}$.
Let $\mathcal{Q}=\{[\rho,\mathbf{m}]|\rho\in[0,\infty),\mathbf{m}\in\mathbb{R}^3\}$ be the natural phase space associated to solutions $[\rho,\mathbf{m}]=[\rho,\rho\mathbf{u}]$.

\subsection{Dissipative measure--valued solutions}

A \emph{dissipative measure-valued (DMV) solution} of the Euler system (1.1) is a parameterized family of probability measures
\begin{align}
\{Y_{t,x}\}_{t\in[0,T],x\in\Omega},\hspace{10pt} (t,x) \mapsto Y_{t,x}\in L^\infty_{weak-(\star)}((0,T)\times\Omega;\mathcal{P}(\mathcal{Q})),
\end{align}
satisfying

\noindent
$\bullet$ the continuity equation
\begin{align}
\int^T_0\int_{\Omega}[\langle Y_{t,x};\rho\rangle\partial_t\varphi+\langle Y_{t,x};\mathbf{m}\rangle\nabla_x\varphi]dxdt
=-\int_{\Omega}\langle Y_{0,x};\rho\rangle\varphi(0)dx,
\end{align}
for all $\varphi\in C^\infty_c([0,T)\times\Omega)$;

\noindent
$\bullet$
the momentum equation
\begin{align}
\int^T_0\int_{\Omega}[\langle Y_{t,x};\mathbf{m}\rangle\partial_t\varphi+\langle Y_{t,x};\frac{\mathbf{m}\otimes\mathbf{m}}{\rho}\rangle:\nabla_x\varphi]dxdt
+\int^T_0\int_{\Omega}\langle Y_{t,x};p(\rho)\rangle\text{div}\varphi dxdt\nonumber\\
+\int^T_0\int_{\Omega}\langle Y_{t,x};\omega\times\mathbf{m}\rangle\varphi dxdt=-\int_{\Omega}\langle Y_{0,x};\mathbf{m}\rangle\varphi(0)dx
-\int^T_0\int_{\Omega}\nabla_x\varphi: d\mu_c,
\end{align}
for all $\varphi\in C^\infty_c([0,T)\times\Omega)$, where  $\mu_c\in\mathcal{M}([0,T]\times\Omega)$
is the so--called momentum concentration measure;

\noindent
$\bullet$
the energy inequality
\begin{align}
\int_{\Omega}[\langle &Y_{\tau,x};\frac{1}{2}\frac{|\mathbf{m}|^2}{\rho}+\big{(}P(\rho)-P'(\overline{\rho})(\rho-\overline{\rho})-P(\overline{\rho})\big{)}\rangle dx+\mathcal{D}(\tau)\nonumber\\
&\leq\int_{\Omega}\langle Y_{0,x};\frac{1}{2}\frac{|\mathbf{m}|^2}{\rho}+\big{(}P(\rho)-P'(\overline{\rho})(\rho-\overline{\rho})-P(\overline{\rho})\big{)}\rangle dx,
\end{align}
for a.a $\tau\in(0,T)$, where
\begin{align}
P(\rho)=\rho\int^\rho_{\overline{\rho}}\frac{p(z)}{z^2}dz,
\end{align}
and $\mathcal{D}$ is a non-negative function $\mathcal{D}\in L^\infty(0,T)$, satisfying the compatibility condition
\begin{align}
\int^\tau_0\int_{\Omega}|\mu_c|dxdt\leq C\int^\tau_0\xi(t)\mathcal{D}(t)dt,\hspace{5pt}\text{for some}\hspace{3pt}\xi\in L^1(0,T).
\end{align}

\begin{remark}
The notion of \emph{(DMV)} solutions can be founded in many  works as it was already mention in the Introduction, see e.g.  \cite{j1,j2,DeStTz,D,e7}. For convenience of readers, we give more details.

Let $L^\infty_{weak\star}((0,T)\times\Omega;P(Q))$ be the space of essentially bounded weakly$^\star-$ measure maps $Y:(0,T)\times\Omega\rightarrow P(Q)$,$(t,x)\mapsto Y_{t,x}$. By virtue of fundamental theorem on Young measures (see \cite{Ba})  there exists a subsequence of $\{\rho_\epsilon,\mathbf m_\epsilon\}_{\epsilon>0}$ and parameterized family of probability measures $\{Y_{t,x}\}_{(t,x)\in(0,T)\times\Omega}$
\begin{eqnarray*}
[(t,x)\mapsto Y_{t,x}]\in L^\infty_{weak\star}((0,T)\times\Omega;P(Q)),
\end{eqnarray*}
such that a.a. $(t,x)\in(0,T)\times\Omega$
\begin{eqnarray*}
\langle Y_{t,x}; G(\rho,\mathbf m)\rangle=\overline{G(\rho,\mathbf m)}(t,x)\hspace{3pt} for\hspace{1pt} any \hspace{3pt} G\in C_c(Q),
\end{eqnarray*}
whenever
\begin{eqnarray*}
G(\rho_\epsilon,\mathbf m_\epsilon)\rightarrow\overline{G(\rho,\mathbf m)}(t,x)\hspace{3pt} weakly\star\hspace{1pt}in \hspace{3pt} L^\infty((0,T)\times\Omega).
\end{eqnarray*}
Moreover, if $G\in C(Q)$ is such that
\begin{eqnarray*}
\int^T_0\int_\Omega |G(\rho_\epsilon,\mathbf m_\epsilon)|dx\leq C,
\end{eqnarray*}
then $G$ is $Y_{t,x}$ integrable for almost all $(t,x)\in(0,T)\times\Omega$ and
\begin{eqnarray*}
[(t,x)\mapsto \langle Y_{t,x}; G(\rho,\mathbf m)\rangle]\in L^1((0,T)\times\Omega),
\end{eqnarray*}
and
\begin{eqnarray*}
G(\rho_\epsilon,\mathbf m_\epsilon)\rightarrow\overline{G(\rho,\mathbf m)}(t,x)\hspace{3pt} weakly\star\hspace{1pt}in \hspace{3pt} \mathcal{M}((0,T)\times\Omega).
\end{eqnarray*}
The difference
\begin{eqnarray*}
\mu_G\equiv\overline{G(\rho,\mathbf m)}-[(t,x)\mapsto \langle Y_{t,x}; G(\rho,\mathbf m)\rangle]\in \mathcal{M}((0,T)\times\Omega),
\end{eqnarray*}
is called concentration defect measure.

For more details, please see \cite{CFKW}.

\end{remark}

\begin{remark}
\smallskip

\begin{itemize}
\item The measure $Y_{0,x}$ plays the role of \emph{initial conditions}.
\item The proof of an existence of \emph{(DMV)} solutions of Euler system was done in the pioneer work  by Neustupa, \cite{Ne}.
Recently, see \cite{n}, the authors proved the local strong solutions of rotating compressible Euler system in $\mathbb{R}^3$. Feireisl et al. \cite{j1,j2} proved the existence of \emph{(DMV)} solutions to the non-rotating full Euler system. As the rotating term does not bring any trouble in the proof of existence, the existence of \emph{(DMV)} solutions to (1.4) can be obtained by analogous methods as in \cite{j2}.

\end{itemize}

\end{remark}

\begin{remark}
We need to define the function
\[
[\rho, \mathbf{m}] \mapsto \frac{ |\mathbf{m}|^2 }{\rho}
\]
on the vacuum set as
\begin{align}
[\rho,\mathbf{m}]\rightarrow\frac{|\mathbf{m}|^2}{\rho}=
\left\{
\begin{array}{llll}  \infty,\hspace{5pt}\text{if} \hspace{3pt} \rho=0 \hspace{3pt} and \hspace{3pt} \mathbf{m}\neq0, \\
\frac{|\mathbf{m}|^2}{\rho},\hspace{5pt}\text{if} \hspace{3pt} \rho>0,\\
0,\hspace{5pt}\text{otherwise}.
\end{array}\right.
\end{align}
Accordingly, it follows from the energy inequality (2.4) that
\begin{align}
{\rm Supp} [Y_{t,x}]\cap\{[\rho,\mathbf{m}]\in\mathcal{Q}|\rho=0,\mathbf{m}\neq0]\}=\emptyset \ \mbox{for a.a.}\ (t,x).
\end{align}
\end{remark}

\subsection{Relative entropy inequality}

Motivated by \cite{e,e1,j2}, we introduce the relative energy functional

\begin{align}
\mathcal{E}(\rho,\mathbf{m}|r, \mathbf{U})=\int_{\Omega}\langle Y_{t,x};\frac{1}{2}\rho|\frac{\mathbf{m}}{\rho}-\mathbf U(t,x)|^2+(P(\rho)-P'(r)(\rho-r)-P(r))\rangle dx,
\end{align}
where $r>0$, $\mathbf U$ are smooth ``test'' functions, $r-\overline{\rho}$, $U$ compactly supported in $\Omega$.

As shown in \cite{j2}, any \emph{(DMV)} solution of (1.1) satisfies the relative entropy inequality
\begin{align}
\mathcal{E}&(\rho,\mathbf{m}|r,\mathbf U)|^{t=\tau}_{t=0}+\mathcal{D(\tau)}\leq
\int^\tau_0\int_{\Omega}\langle Y_{t,x};(\partial_t \mathbf U+\frac{\mathbf{m}}{\rho}\nabla_xU)(\rho \mathbf U-\mathbf{m})\rangle dxdt\nonumber\\
&+\int^\tau_0\int_{\Omega}\langle Y_{t,x};(r-\rho)\partial_tP'(r)
+(r \mathbf U-\mathbf{m})\nabla_xP'(r)\rangle dxdt+\int^\tau_0\int_{\Omega}\langle Y_{t,x};\omega\times\frac{\mathbf{m}}{\rho}\rangle(\rho\mathbf{U}-\mathbf{m})dxdt\nonumber\\
&-\int^\tau_0\int_{\Omega}\langle Y_{t,x};p(\rho)-p(r)\rangle\text{div} \mathbf Udxdt
+\int^\tau_0\int_{\Omega}\nabla_x \mathbf U:d\mu_c.
\end{align}
for a.a. $\tau \in [0,T]$,
and any $r, \mathbf U\in$$C^1([0,T]\times\Omega)$, $r-\overline{\rho}$, $\mathbf U$ compactly supported in $\Omega$.
\section{Main result}
Before stating our main result, we introduce some notations and  collect several mostly technical hypotheses and known facts concerning the limit system. $x=(x_h,x_3)$ with $x_h\in\mathbb{R}^2$ denoting its horizontal component. For a vector field $\mathbf b=[b_1,b_2,b_3]$, we introduce the horizontal component $\mathbf b_h=[b_1,b_2]$ writing $\mathbf b=[\mathbf b_h,b_3]$. Similarly, we use the symbols $\nabla_h$, $\text{div}_h$ to denote the differential operators acting on the horizontal variables. The following assumptions and results will be used in the proof.
\subsection{Pressure}

We suppose the pressure $p$ is a continuously differentiable function of the density such that for some $\gamma>1$,
\begin{equation} \label{pres}
p\in C^1[0,\infty)\cap C^\infty(0,\infty),\hspace{3pt}p(0)=0,\hspace{3pt}p'(\rho)>0\hspace{3pt} \text{for all}\hspace{3pt}\rho>0,
\lim_{\rho\rightarrow\infty}\frac{p'(\rho)}{\rho^{\gamma-1}}=p_\infty >0.
\end{equation}

\begin{remark}
{Similarly to \cite{e6},} we deduce that
\begin{align*}
&p(\rho)-p'(r)(\rho-r)-p(r)\hspace{3pt}\text{is dominated by}\hspace{3pt}P(\rho)-P'(r)(\rho-r)-P(r),\ \mbox{specifically,}\\
&|\rho-r|^2\leq c(\delta)(P(\rho)-P'(r)(\rho-r)-P(r))\hspace{8pt}\text{when}\hspace{3pt} 0<\delta\leq\rho,\hspace{3pt}r\leq\frac{1}{\delta},\hspace{3pt}\delta>0,\\
&{1+|\rho-r|+P(\rho)}\leq c(\delta)(P(\rho)-P'(r)(\rho-r)-P(r))\hspace{8pt}\text{if}\hspace{3pt} 0<2\delta<r<\frac{1}{2\delta},\\
&\hspace{245pt}\rho\in[0,\delta)\cup(\frac{1}{\delta},\infty),\hspace{3pt}\delta>0.
\end{align*}
\end{remark}

\subsection{Quasi-geophysical equation}

The expected limit problem reads
\begin{align}
&\omega\times\mathbf{v}+\frac{p'(\overline{\rho})}{\overline{\rho}}\nabla_xq=0,\hspace{5pt}\mathbf{v}=[\mathbf{v_h}(x_h),0],\hspace{2pt}q=q(x_h),\\
&\partial_t(\Delta_hq-\frac{1}{p'(\overline{\rho})}q)+\mathbf{v}_h\cdot\nabla_h(\Delta_hq)=0.
\end{align}
supplement with the initial condition
\begin{align*}
q|_{t=0}=q_0.
\end{align*}
As shown by Oliver \cite{o}, the problem $(3.2)-(3.3)$ possesses a unique classical solution
\begin{align}
q\in C([0,T];W^{m,2}(\mathbb{R}^2))\cap C^1([0,T];W^{m-1,2}(\mathbb{R}^2)),\hspace{8pt}m\geq4,
\end{align}
for any initial solution
\begin{align}
q_0\in W^{m,2}(\mathbb{R}^2).
\end{align}

\subsection{Ill prepared initial--data}
\label{ILP}
The \emph{ill--prepared} initial data for the scaled system \eqref{1.4} take the form
\begin{eqnarray}
\rho_\epsilon (0, \cdot) = \rho_{0, \epsilon} =\overline{\rho} + \epsilon s_{0, \epsilon},\
\mathbf{u}_\epsilon (0, \cdot) = \mathbf{u}_{0, \epsilon},
\end{eqnarray}
where
\begin{eqnarray}
&s_{0,\epsilon} \rightarrow s_0\hspace{3pt}\text{in}\hspace{3pt}W^{k,2}(\Omega)\cap W^{k,1}(\Omega),\hspace{8pt}
\mathbf{u}_{0,\epsilon}\rightarrow \mathbf{u}_0\hspace{3pt}\text{in}\hspace{3pt}W^{k,2}(\Omega)\cap W^{k,1}(\Omega),\hspace{5pt}(k>3),\\
&\mathbf{u}_0 = \mathbf{v}_0 + \nabla_x \Phi_0{\color{blue}.}
\end{eqnarray}

\subsection{Singular limit -- main result}

For simplicity, we assume $\overline{\rho}=p'(\overline{\rho})=P''(\overline{\rho})=1$. Now, we are ready to state our main result.
\begin{theorem}
Let $\{Y_{t,x}^\epsilon\}_{(t,x)\in[0,T]\times\Omega}$ be a family of (DMV) solutions to the scaled Euler system \eqref{1.4} satisfying the compatibility condition (2.6)
with a function $\xi$ independent of $\epsilon$.
Let the initial data $\{ Y_{0,x}^\epsilon \}_{x \in \Omega}$ be ill-prepared, namely
\begin{align*}
\int_{\Omega}\langle Y_{0,x}^\epsilon;\frac{1}{2}\rho|\frac{\mathbf{m}}{\rho}-\mathbf{u}_{0, \epsilon} (x)|^2+\frac{1}{\epsilon^2}(P(\rho)-P'(\rho_{0, \epsilon})(\rho- \rho_{0, \epsilon})-P(\rho_{0, \epsilon}))\rangle dx\rightarrow0,
\end{align*}
where $\rho_{0, \epsilon}$, $\mathbf{u}_{0, \epsilon}$ are ill prepared data introduced in Section \ref{ILP}.

Then
\begin{align*}
&\mathcal{D}^\epsilon\rightarrow0\hspace{5pt}\text{in}\hspace{3pt}L^\infty(0,T),\\
&Y_{t,x}^\epsilon\rightarrow\delta_{[q,\mathbf{v}]}\hspace{5pt}\text{in}\hspace{3pt}L^p(0,T;L^1_{{\rm loc}}(\Omega; \mathcal{M}^+(\mathcal{Q})_{weak-(\star)}))\ \mbox{for any finite}\ p \geq 1,
\end{align*}
where $q$ and $\mathbf{v}$ is the unique solution of problem (3.2)-(3.3) starting from the initial data $q_0$ and where $q_0\in W^{k+1,2}(\mathbb{R}^2)\cap W^{k+1,1}(\mathbb{R}^2)$ is the unique solution of the elliptic problem
\begin{align}
-\Delta_hq_0+q_0=\int^1_0\text{curl}_h[\mathbf{u}_0]_hdx_3+\int^1_0s_0dx_3.
\end{align}

\end{theorem}

The rest of the paper is devoted to the proof of Theorem 3.1.

\section{Energy bounds}

We start by deriving uniform bounds on solutions to \eqref{1.4} independent of $\epsilon$.  Similarly to \cite{e6}, we introduce the decomposition
\begin{align*}
h (\rho, \mathbf{m} ) =[h]_{ess} (\rho, \mathbf{m}) +[h]_{res} (\rho, \mathbf{m})   ,\hspace{10pt}[h]_{ess}=\psi(\rho )h(\rho, \mathbf{m}),\hspace{5pt}[h]_{res}=(1-\psi(\rho))h(\rho, \mathbf{m}),
\end{align*}
where
\begin{align*}
\psi \in C_c^\infty(0,\infty),\hspace{5pt}0\leq\psi(\rho)\leq1,\hspace{5pt}\psi(\rho)=1\hspace{3pt}\text{on an open interval containing}\hspace{3pt}\overline{\rho}=1.
\end{align*}

As the initial data are ill--prepared, the expression on the right--hand side of the energy inequality (2.4) remains bounded uniformly for $\epsilon \to 0$.
Consequently, we deduce the following bound:
\begin{align}
ess\sup_{t\in(0,T)}\int_{\Omega}\langle Y_{t,x}^\epsilon;\frac{1}{2}\frac{|\mathbf{m}_\epsilon|^2}{\rho}+\frac{1}{\epsilon^2}(P(\rho_\epsilon)-P'(1)(\rho_\epsilon-1)-P(1))\rangle dx\leq C.
\end{align}

Thus, {exactly as in \cite{e6},} we use the structural properties of the function $p$ to deduce
\begin{align}
&ess\sup_{t\in(0,T)}\int_{\Omega}\langle Y_{t,x}^\epsilon;|[\frac{\rho_\epsilon-1}{\epsilon}]_{ess}|^2\rangle+
\langle Y_{t,x}^\epsilon;[\frac{P(\rho_\epsilon)+1}{\epsilon^2}]_{ess}\rangle dx\leq C;\nonumber\\
&(t,x) \mapsto\langle Y_{t,x}^\epsilon;\mathbf{m}_\epsilon\rangle\hspace{5pt}\text{bounded in}\hspace{3pt} L^\infty(0,T; L^2(\Omega) +  L^{\frac{2\gamma}{\gamma+1}}(\Omega));\nonumber\\
&(t,x) \mapsto \langle Y_{t,x}^\epsilon;[\frac{\rho_\epsilon-1}{\epsilon}]_{ess}\rangle\hspace{5pt}\text{bounded in}\hspace{3pt} L^\infty(0,T; L^2(\Omega));\nonumber\\
&(t,x) \mapsto \epsilon^{-\frac{2}{\gamma}}\langle Y_{t,x}^\epsilon;[\rho_\epsilon]_{res}\rangle\hspace{5pt}\text{bounded in}\hspace{3pt} L^\infty(0,T; L^\gamma(\Omega)).
\end{align}

Using the same argument in \cite{e5}, there exist functions $\rho^{(1)}\in L^\infty(0,T;L^2(\Omega))$ and $\mathbf m\in L^\infty(0,T;L^q(\Omega))$ for some $q>1$ and a subsequence such that
\begin{align*}
&\langle Y_{t,x}^\epsilon;\mathbf{m}_\epsilon\rangle\rightarrow\mathbf m\hspace{5pt}\text{weakly in}\hspace{3pt} L^\infty(0,T; L^q(\Omega);\\
&\langle Y_{t,x}^\epsilon;\frac{\rho_\epsilon-1}{\epsilon}\rangle\rightarrow \rho^{(1)}\hspace{5pt}\text{weakly in}\hspace{3pt} L^\infty(0,T; L^2(\Omega)).
\end{align*}

Recalling (2.2) and (2.3), we deduce
\begin{align*}
\int^T_0\int_\Omega\mathbf m\cdot\nabla_x\varphi dxdt=0,\hspace{5pt}\int^T_0\int_\Omega[(\omega\times\mathbf m)\cdot\varphi+\rho^{(1)}\text{div}\varphi]dxdt=0,
\end{align*}
for $\varphi\in C^1([0,T]\times\Omega)$. In other words,
\begin{align}
\text{div}_x\mathbf m=0,\hspace{5pt}\omega\times\mathbf m+\nabla_x\rho^{(1)}=0,
\end{align}
in the sense of distribution.

It is easy to check that
\begin{align*}
\rho^{(1)}=\rho^{(1)}(x_h),\hspace{3pt}\mathbf m=(\mathbf m_h,0),\hspace{3pt}
\text{div}_x\mathbf m=\text{div}_h\mathbf m_h=0.
\end{align*}

Moreover, the detail of derivation of (3.9) can be seen in \cite{e2,e4,no}.
\section{Acoustic waves}
It is well-known that ill-prepared data give rise to rapidly oscillating acoustic waves. Similarly to \cite{e4}, the relevant acoustic equation reads
\begin{eqnarray}
\left\{
\begin{array}{llll}  \epsilon\partial_{t}s_\epsilon+\text{div}(\nabla_x\Phi_\epsilon)=0, \\
\epsilon\partial_t\nabla_x\Phi_\epsilon+\omega\times\nabla_x\Phi_\epsilon+\nabla_xs_\epsilon=0,
\end{array}\right.
\end{eqnarray}
supplemented with the initial data
\[
s_\epsilon (0, \cdot) = s_{0},\ \nabla_x \Phi_\epsilon (0, \cdot) = \nabla_x \Phi_0,
\]
where $s_0$, $\nabla_x \Phi_0$ have been introduced in Section \ref{ILP}.

As a matter of fact, the initial data must be smoothed and cut-off via suitable regularization operators, namely
\begin{align}
&s_\epsilon(0,\cdot)=s_{0, \delta}=[s_{0}]_{\delta};\hspace{5pt}
\nabla_x \Phi_\epsilon(0,\cdot)= \nabla_x \Phi_{0, \delta}= \nabla_x [\Phi_{0}]_{\delta},\nonumber
\end{align}
where $[\cdot]_\delta$ denotes the regularization introduced in \cite{e4}.

Denoting the corresponding solutions $s_{\epsilon, \delta}$, $\Phi_{\epsilon, \delta}$ we report the following energy and dispersive estimates proved in \cite[Section 6]{e4}:
\begin{align}
\sup_{t\in[0,T]}[\|\Phi_{\epsilon , \delta} (t,\cdot)\|_{W^{m,2}}+\|s_{\epsilon, \delta} (t,\cdot)\|_{W^{m,2}}]=[\|\nabla_x\Phi_{0,\delta}\|_{L^{2}}+\|s_{0,\delta}\|_{L^{2}}],
\end{align}
and
\begin{align}
{\int^T_0[\|\Phi_{\epsilon, \delta} (t,\cdot)\|_{W^{m,\infty}}+\|s_{\epsilon, \delta} (t,\cdot)\|_{W^{m,\infty}}]\leq
\omega(\epsilon, m , \delta)[\|\nabla_x\Phi_{0,\delta}\|_{L^{2}}+\|r_{0,\delta}\|_{L^{2}}],}
\end{align}
where $\omega(\epsilon, m, \delta )\rightarrow0$ as $\epsilon\rightarrow0$ for any fixed $m \geq 0$ and $\delta > 0$. More details about Strichartz estimates and acoustic waves, readers can refer to \cite{s,sm}.
\section{Convergence}

The proof of convergence is based on the ansatz
\begin{equation}
r_\epsilon = 1 + \epsilon(q+s_{\epsilon, \delta}),\ \mathbf{U}_\epsilon = \mathbf{v} + \nabla_x \Phi_{\epsilon, \delta},
\end{equation}
in the relative energy inequality (2.10). The $[s_{\epsilon, \delta},\nabla_x\Phi_{\epsilon, \delta}]$ are solutions of the acoustic system (5.1), and $[q,\mathbf v]$ is solution of the target problem
\begin{align}
&\omega\times\mathbf{v}+\nabla_xq=0,\nonumber\\
&\partial_t(\Delta_hq-q)+\nabla_h^{\perp}q\cdot\nabla_h(\Delta_hq)=0.
\end{align}
In addition, to avoid technicalities, we shall assume that $s_0$ and $\Phi_0$ are sufficiently regular so that
the $\delta-$regularization is not needed in (5.1--5.3). Accordingly, we have $s_{\epsilon, \delta} = s_\epsilon$, $\Phi_{\epsilon, \delta} = \Phi_{\epsilon}$. The general case may be handled as in \cite{e4}.

First note that the relative energy for the scaled system reads
\begin{align}
\mathcal{E}_\epsilon (\rho_\epsilon,\mathbf{m}_\epsilon|r_\epsilon, \mathbf U_\epsilon)=\int_{\Omega}\langle Y_{t,x};\frac{1}{2}\rho_\epsilon|\frac{\mathbf{m}_\epsilon}{\rho_\epsilon}-\mathbf U_\epsilon|^2+\frac{1}{\epsilon^2}(P(\rho_\epsilon)-P'(r_\epsilon)(\rho_\epsilon-r_\epsilon)-P(r_\epsilon))\rangle dx,
\end{align}
with the corresponding relative energy inequality:
\begin{align}
\mathcal{E}_\epsilon &(\rho_\epsilon,\mathbf{m}_\epsilon|r_\epsilon, \mathbf U_\epsilon)|^{t=\tau}_{t=0}+\mathcal{D}^{\epsilon}(\tau)\leq
\int^\tau_0\int_{\Omega}\langle Y^\epsilon_{t,x};\rho_\epsilon \mathbf U_\epsilon-\mathbf{m}_\epsilon\rangle(\partial_t \mathbf U_\epsilon+\frac{\mathbf{m}_\epsilon}{\rho_\epsilon}\nabla_x \mathbf U_\epsilon)dxdt\nonumber\\
&+\frac{1}{\epsilon^2}\int^\tau_0\int_{\Omega}[\langle Y^\epsilon_{t,x};r_\epsilon-\rho_\epsilon\rangle\partial_tP'(r_\epsilon)
+\langle Y^\epsilon_{t,x};r_\epsilon \mathbf U_\epsilon-\mathbf{m}_\epsilon\rangle\nabla_xP'(r_\epsilon) ]dxdt\nonumber\\
&+ \frac{1}{\epsilon}\int^\tau_0\int_{\Omega}\langle Y_{t,x};\omega\times\frac{\mathbf{m}_\epsilon}{\rho_\epsilon}\rangle(\rho_\epsilon\mathbf{U}_\epsilon-\mathbf{m}_\epsilon)dxdt-\frac{1}{\epsilon^2} \int^\tau_0\int_{\Omega}\langle Y_{t,x};p(\rho_\epsilon)-p(r_\epsilon)\rangle\text{div} \mathbf U_\epsilon dxdt\nonumber\\
&+\int^\tau_0\int_{\Omega}\nabla_x \mathbf U_\epsilon :d\mu_c.
\end{align}

Our goal is to show that, with the ansatz (6.1), the relative energy $\mathcal{E}_\epsilon (\rho_\epsilon,\mathbf{m}_\epsilon|r_\epsilon, \mathbf U_\epsilon)$ tends to
zero for $\epsilon \to 0$ uniformly in $t \in [0,T]$. In view of the dispersive estimates $(5.2)-(5.3)$, this will yield the conclusion claimed in Theorem 3.1.
To this end, we use a Gronwall type argument
showing that all integrals in the right-hand side of (6.4) are either small or can be absorbed by the left-hand side as $\epsilon\rightarrow0$.
This programme will be carried over by means of several steps.

\subsection{Step 1} {First, we compute
\begin{align*}
&\int^\tau_0\int_{\Omega}[\langle Y^\epsilon_{t,x};r_\epsilon-\rho_\epsilon\rangle\partial_tP'(r_\epsilon)
+\langle Y^\epsilon_{t,x};r_\epsilon \mathbf U_\epsilon-\mathbf{m}_\epsilon\rangle\nabla_xP'(r_\epsilon)-\langle Y^\epsilon_{t,x};p(\rho_\epsilon)-p(r_\epsilon)\rangle\text{div} \mathbf U_\epsilon]dxdt\\
&=\int^\tau_0\int_{\Omega}[\langle Y^\epsilon_{t,x};p(r_\epsilon)-p'(r_\epsilon)(r_\epsilon-\rho_\epsilon)-p(\rho_\epsilon)\rangle\text{div} \mathbf U_\epsilon
+\langle Y^\epsilon_{t,x};r_\epsilon-\rho_\epsilon\rangle\partial_tP'(r_\epsilon)\\
&\hspace{8pt}+\langle Y^\epsilon_{t,x};(r_\epsilon-\rho_\epsilon)p'(r_\epsilon)\rangle\text{div}\mathbf U_\epsilon
+\langle Y^\epsilon_{t,x};(r_\epsilon-\rho_\epsilon)\nabla_xP'(r_\epsilon)\rangle \mathbf U_\epsilon+\langle Y^\epsilon_{t,x};(\rho_\epsilon\mathbf U_\epsilon-\mathbf m_\epsilon)\nabla_xP'(r_\epsilon)\rangle]dxdt\\
&=\int^\tau_0\int_{\Omega}[\langle Y^\epsilon_{t,x};p(r_\epsilon)-p'(r_\epsilon)(r_\epsilon-\rho_\epsilon)-p(\rho_\epsilon)\rangle\text{div}  \mathbf U_\epsilon
+\langle Y^\epsilon_{t,x};\partial_tr_\epsilon+\text{div}_x(r_\epsilon \mathbf U_\epsilon)\rangle(r_\epsilon-\rho_\epsilon)P''(r_\epsilon)\\
&\hspace{8pt}+\langle Y^\epsilon_{t,x};(\rho_\epsilon\mathbf U_\epsilon-\mathbf m_\epsilon)\nabla_xP'(r_\epsilon)\rangle]dxdt.
\end{align*}

Note that, in view of (6.2),}
\begin{align*}
\partial_t r_\epsilon +\text{div}_x(r_\epsilon \mathbf U_\epsilon)&=\epsilon\partial_tq+\partial_ts_\epsilon+\text{div}(r_\epsilon(\mathbf v+\nabla_x\Phi_\epsilon))\\
&=\epsilon\partial_tq+\epsilon\text{div}((q+s_\epsilon)\mathbf{U}_\epsilon)).
\end{align*}

{Next, by virtue of (5.1) and (6.1),}
\begin{align*}
\nabla_xP'(r_\epsilon)(\rho_\epsilon\mathbf U_\epsilon-\mathbf{m}_\epsilon)&=\nabla_x(P'(r_\epsilon)-P''(1)(r_\epsilon-1)-P'(1))(\rho_\epsilon\mathbf U_\epsilon-\mathbf{m}_\epsilon)+\epsilon\nabla_xq\cdot(\rho_\epsilon\mathbf U_\epsilon-\mathbf{m}_\epsilon)\\
&\hspace{20pt}+\epsilon\nabla_xs_\epsilon\cdot(\rho_\epsilon\mathbf U_\epsilon-\mathbf{m}_\epsilon)\\
&=\nabla_x(P'(r_\epsilon)-P''(1)(r_\epsilon-1)-P'(1))(\rho_\epsilon\mathbf U_\epsilon-\mathbf{m}_\epsilon)+\epsilon\nabla_xq\cdot(\rho_\epsilon\mathbf U_\epsilon-\mathbf{m}_\epsilon)\\
&\hspace{20pt}-\epsilon^2(\rho_\epsilon\mathbf U_\epsilon-\mathbf{m}_\epsilon)\cdot\partial_t\nabla\Phi_\epsilon-\epsilon(\rho_\epsilon\mathbf U_\epsilon-\mathbf{m}_\epsilon)(\omega\times\nabla\Phi_\epsilon).
\end{align*}

Furthermore, by virtue of the compatibility condition (2.6), we can control the concentration measure,
\begin{align*}
{\int^\tau_0\int_{\Omega}\nabla_x \mathbf U:d\mu_c}\leq\|\nabla_x \mathbf U\|_{L^\infty}\int^\tau_0\xi(t)\mathcal{D}^\epsilon(t)dt.
\end{align*}

{Finally,
as the hypotheses about the ill-prepared initial data,} we have
\[
\mathcal{E}_\epsilon (\rho_\epsilon,\mathbf{m}_\epsilon|r_\epsilon, \mathbf U_\epsilon)(0) \to 0 \ \mbox{as}\ \epsilon \to 0.
\]

Thus we may conclude that
\begin{align*}
\mathcal{E}_\epsilon &(\rho_\epsilon,\mathbf{m}_\epsilon|r_\epsilon, \mathbf U_\epsilon)(\tau) +\mathcal{D}^{\epsilon}(\tau)\leq
\int^\tau_0\int_{\Omega}\langle Y^\epsilon_{t,x};(\partial_t \mathbf{v}+\frac{\mathbf{m}_\epsilon}{\rho_\epsilon}\nabla_x \mathbf U_\epsilon)(\rho_\epsilon \mathbf U_\epsilon-\mathbf{m}_\epsilon)\rangle dxdt\\
&+\frac{1}{\epsilon}\int^\tau_0\int_{\Omega}\langle Y^\epsilon_{t,x};\omega\times\mathbf{v}\rangle(\rho_\epsilon \mathbf U_\epsilon-\mathbf{m}_\epsilon)dxdt\\
&+\frac{1}{\epsilon^2}\int^\tau_0\int_{\Omega}\langle Y^\epsilon_{t,x};\nabla_x\big{(}P'(r_\epsilon)-P''(1)(r_\epsilon-1)-P'(1)\big{)}\rangle(\rho_\epsilon \mathbf U_\epsilon-\mathbf{m}_\epsilon)dxdt\\
&-\frac{1}{\epsilon^2}\int^\tau_0\int_{\Omega}\langle Y^\epsilon_{t,x};p(\rho_\epsilon)-p(r_\epsilon)-p'(r_\epsilon)(\rho_\epsilon-r_\epsilon)\rangle\text{div} \mathbf U_\epsilon dxdt\\
&+\frac{1}{\epsilon}\int^\tau_0\int_{\Omega}\langle Y^\epsilon_{t,x};\partial_tq+\text{div}((q+s_\epsilon)\mathbf{U}_\epsilon))\rangle(r_\epsilon-\rho_\epsilon)P''(r_\epsilon)dxdt\\
&+\frac{1}{\epsilon}\int^\tau_0\int_{\Omega}\langle Y^\epsilon_{t,x};\nabla_xq\rangle(\rho_\epsilon \mathbf U_\epsilon-\mathbf{m}_\epsilon)dxdt+c\int^\tau_0\xi(t)\mathcal{D}^\epsilon(t)dt + \omega(\epsilon),
\end{align*}
where $\omega(\epsilon)$ denotes a generic quantity satisfying
\[
\omega(\epsilon) \to 0 \ \mbox{in}\ L^1(0,T)\ \mbox{as}\ \epsilon \to 0.
\]
Using (6.2), we get the following conclusion:
\begin{align*}
\mathcal{E}_\epsilon &(\rho_\epsilon,\mathbf{m}_\epsilon|r_\epsilon, \mathbf U_\epsilon)(\tau) +\mathcal{D}^{\epsilon}(\tau)\leq
\int^\tau_0\int_{\Omega}\langle Y^\epsilon_{t,x};(\partial_t \mathbf{v}+\frac{\mathbf{m}_\epsilon}{\rho_\epsilon}\nabla_x \mathbf U_\epsilon)(\rho_\epsilon \mathbf U_\epsilon-\mathbf{m}_\epsilon)\rangle dxdt\\
&+\frac{1}{\epsilon^2}\int^\tau_0\int_{\Omega}\langle Y^\epsilon_{t,x};\nabla_x\big{(}P'(r_\epsilon)-P''(1)(r_\epsilon-1)-P'(1)\big{)}\rangle(\rho_\epsilon \mathbf U_\epsilon-\mathbf{m}_\epsilon)dxdt\\
&-\frac{1}{\epsilon^2}\int^\tau_0\int_{\Omega}\langle Y^\epsilon_{t,x};p(\rho_\epsilon)-p(r_\epsilon)-p'(r_\epsilon)(\rho_\epsilon-r_\epsilon)\rangle\text{div} \mathbf U_\epsilon dxdt\\
&+\frac{1}{\epsilon}\int^\tau_0\int_{\Omega}\langle Y^\epsilon_{t,x};\partial_tq+\text{div}((q+s_\epsilon)\mathbf{U}_\epsilon))\rangle(r_\epsilon-\rho_\epsilon)P''(r_\epsilon)dxdt+c\int^\tau_0\xi(t)\mathcal{D}^\epsilon(t)dt + \omega(\epsilon),
\end{align*}

\subsection{Step 2}
We write
\begin{align*}
\int^\tau_0\int_{\Omega}&[\langle \mathbf Y^\epsilon_{t,x};\rho_\epsilon U_\epsilon-\mathbf{m}_\epsilon\rangle(\partial_t \mathbf{v}+\frac{\mathbf{m}_\epsilon}{\rho_\epsilon}\nabla_x \mathbf U_\epsilon)]dxdt\\
&=\int^\tau_0\int_{\Omega}\langle Y^\epsilon_{t,x};\rho_\epsilon \mathbf U_\epsilon-\mathbf{m}_\epsilon\rangle(\partial_t\mathbf{v}+\mathbf{v}\cdot\nabla_x\mathbf{v})dxdt\\
&\hspace{8pt}+\int^\tau_0\int_{\Omega}\langle Y^\epsilon_{t,x};\rho_\epsilon \mathbf U_\epsilon-\mathbf{m}_\epsilon\rangle(\mathbf{v}\cdot\nabla_x\nabla_x\Phi_\epsilon+\nabla_x\Phi_\epsilon\nabla_x \mathbf U_\epsilon)dxdt\\
&\hspace{8pt}+\int^\tau_0\int_{\Omega}\langle Y^\epsilon_{t,x};\rho_\epsilon \mathbf U_\epsilon-\mathbf{m}_\epsilon\rangle(\frac{\mathbf{m}_\epsilon}{\rho_\epsilon}-\mathbf U_\epsilon)\nabla_x \mathbf U_\epsilon dxdt\\
&=I_1+I_2+I_3.
\end{align*}

{Using the uniform bounds (4.2)}, we can split the functions in $I_2$ into their essential and residual parts obtaining
\[
\begin{split}
&\left| \int_{\mathbb{R}^2}\langle Y^\epsilon_{t,x};\rho_\epsilon \mathbf U_\epsilon-\mathbf{m}_\epsilon\rangle(\mathbf{v}\cdot\nabla_x\nabla_x\Phi_\epsilon+\nabla_x\Phi_\epsilon\nabla_x \mathbf U_\epsilon) \ dx \right| \\
&\leq \| \nabla_x \Phi_\epsilon \|^2_{W^{1,\infty}} \left( \| \mathbf{v} \|_{W^{3,2}} + \| \nabla_x \mathbf{U}_\epsilon \|_{W^{3,2}} \right)^2 +
c  \mathcal{E}_\epsilon (\rho_\epsilon,\mathbf{m}_\epsilon|r_\epsilon, \mathbf U_\epsilon),
\end{split}
\]
where the first term on the right--hand side can be controlled by means of the dispersive estimate (5.2) and (5.3).

Summing up the previous observations, we may infer that the relative energy inequality with the ansatz (6.1) reduces to
\begin{align*}
\mathcal{E}_\epsilon &(\rho_\epsilon,\mathbf{m}_\epsilon|r_\epsilon, \mathbf U_\epsilon)(\tau) +\mathcal{D}^{\epsilon}(\tau)\leq\int^\tau_0\int_{\Omega}\langle Y^\epsilon_{t,x};\rho_\epsilon \mathbf U_\epsilon-\mathbf{m}_\epsilon\rangle(\partial_t\mathbf{v}+\mathbf{v}\cdot\nabla_x\mathbf{v})dxdt
\\
&+\frac{1}{\epsilon^2}\int^\tau_0\int_{\Omega}\langle Y^\epsilon_{t,x};\nabla_x\big{(}P'(r_\epsilon)-P''(1)(r_\epsilon-1)-P'(1)\big{)}\rangle(\rho_\epsilon \mathbf U_\epsilon-\mathbf{m}_\epsilon)dxdt\\
&-\frac{1}{\epsilon^2}\int^\tau_0\int_{\Omega}\langle Y^\epsilon_{t,x};p(\rho_\epsilon)-p(r_\epsilon)-p'(r_\epsilon)(\rho_\epsilon-r_\epsilon)\rangle\text{div} \mathbf U_\epsilon dxdt\\
&-\frac{1}{\epsilon}\int^\tau_0\int_{\Omega}\langle Y^\epsilon_{t,x};(\rho_\epsilon-r_\epsilon)P''(r_\epsilon)\rangle(\partial_tq+\text{div}((q+s_\epsilon)\mathbf{U}_\epsilon))dxdt\\
&+C\int^\tau_0\mathcal{E}_\epsilon (\rho_\epsilon,\mathbf{m}_\epsilon|r_\epsilon, \mathbf U_\epsilon)dt+C\int^\tau_0\xi(t)\mathcal{D}^\epsilon(t)dt + \omega(\epsilon).
\end{align*}

\subsection{Step 3}
Now, we will deal with pressure term and corresponding term. First, using direct calculation, the Taylor formula and dispersive estimates (5.2-5.3), we deduce that
\begin{align*}
\frac{1}{\epsilon^2}|\nabla_x\big{(}P'(r_\epsilon)-&P'(1)-P''(1)(r_\epsilon-1)\big{)}|\\
&=\frac{1}{\epsilon}|\big{(}P''(r_\epsilon)-P''(1)\big{)}\nabla_x(q+s_\epsilon)|\\
&\rightarrow P'''(1)q\nabla_xq \hspace{5pt}\text{as}\hspace{3pt}\epsilon\rightarrow0.
\end{align*}

Therefore, combining the previous energy bounds and convergence, we get
\begin{align*}
\frac{1}{\epsilon^2}\int^\tau_0\langle Y^\epsilon_{t,x};&\nabla_x\big{(}P'(r_\epsilon)-P'(1)-P''(1)(r_\epsilon-1)\big{)}\rangle(\rho_\epsilon \mathbf U_\epsilon-\mathbf{m}_\epsilon)dt\rightarrow 0\hspace{5pt}\text{as}\hspace{3pt}\epsilon\rightarrow0.
\end{align*}

The remaining pressure term is
\begin{align*}
|\frac{1}{\epsilon^2}&\int^\tau_0\langle Y^\epsilon_{t,x};p(\rho_\epsilon)-p(r_\epsilon)-p'(r_\epsilon)(\rho_\epsilon-r_\epsilon)\rangle\text{div} \mathbf U_\epsilon dt|\\
&=|\frac{1}{\epsilon^2}\int^\tau_0\langle Y^\epsilon_{t,x};p(\rho_\epsilon)-p(r_\epsilon)-p'(r_\epsilon)(\rho_\epsilon-r_\epsilon)\rangle(\text{div}\mathbf{v}+\Delta\Phi_\epsilon)dt|\\
&\leq c |\frac{1}{\epsilon^2}\int^\tau_0\langle Y^\epsilon_{t,x};P(\rho_\epsilon)-P(r_\epsilon)-P'(r_\epsilon)(\rho_\epsilon-r_\epsilon)\rangle(\text{div}\mathbf{v}+\Delta\Phi_\epsilon)dt|\\
&\leq C\int^\tau_0\mathcal{E}_\epsilon (\rho_\epsilon,\mathbf{m}_\epsilon|r_\epsilon, \mathbf U_\epsilon)dt,
\end{align*}
where we have used the previous dispersive estimates (5.2) and (5.3). Thus, we can conclude that
\begin{align*}
\mathcal{E}_\epsilon &(\rho_\epsilon,\mathbf{m}_\epsilon|r_\epsilon, \mathbf U_\epsilon)(\tau) +\mathcal{D}^{\epsilon}(\tau)\leq \omega(\epsilon) +C\int^\tau_0\mathcal{E}(\rho_\epsilon,\mathbf{m}_\epsilon|r_\epsilon,U_\epsilon)dt+c\int^\tau_0\xi(t)\mathcal{D}^\epsilon(t)dt\\
&+\int^\tau_0\int_{\Omega}\langle Y^\epsilon_{t,x};\rho_\epsilon \mathbf U_\epsilon-\mathbf{m}_\epsilon\rangle(\partial_t\mathbf{v}+\mathbf{v}\cdot\nabla_x\mathbf{v})dxdt\\
&-\frac{1}{\epsilon}\int^\tau_0\int_{\Omega}\langle Y^\epsilon_{t,x};(\rho_\epsilon-r_\epsilon)P''(r_\epsilon)\rangle(\partial_tq+\text{div}((q+s_\epsilon)\mathbf{U}_\epsilon))dxdt.
\end{align*}
\subsection{Step 4} Finally, we deal with the remaining pressure terms.
Similar to the previous analysis, we obtain
\begin{align*}
\mathcal{E}_\epsilon &(\rho_\epsilon,\mathbf{m}_\epsilon|r_\epsilon, \mathbf U_\epsilon)(\tau) +\mathcal{D}^{\epsilon}(\tau)\leq \omega(\epsilon) +C\int^\tau_0\mathcal{E}(\rho_\epsilon,\mathbf{m}_\epsilon|r_\epsilon,U_\epsilon)dt+c\int^\tau_0\xi(t)\mathcal{D}^\epsilon(t)dt\\
&+\int^\tau_0\int_{\Omega}\langle Y^\epsilon_{t,x};\mathbf v-\mathbf{m}\rangle(\partial_t\mathbf{v}+\mathbf{v}\cdot\nabla_x\mathbf{v})dxdt
+\int^\tau_0\int_{\Omega}\langle Y^\epsilon_{t,x};q-\rho^{(1)}\rangle(\partial_tq+\text{div}(q\mathbf{v}))dxdt,
\end{align*}
where
\begin{align*}
&\int^\tau_0\int_{\Omega}\langle Y^\epsilon_{t,x};\mathbf v-\mathbf{m}\rangle(\partial_t\mathbf{v}+\mathbf{v}\cdot\nabla_x\mathbf{v})dxdt
+\int^\tau_0\int_{\Omega}\langle Y^\epsilon_{t,x};q-\rho^{(1)}\rangle(\partial_tq+\text{div}(q\mathbf{v}))dxdt\\
&=\frac{1}{2}\int^\tau_0\int_{\Omega}\langle Y^\epsilon_{t,x};\partial_t|\mathbf v|^2+\partial_t|q|^2\rangle dxdt
-\int^\tau_0\int_{\Omega}\langle Y^\epsilon_{t,x};\partial_t\mathbf v\cdot\mathbf m+\partial_tq\cdot\rho^{(1)}\rangle dxdt\\
&\hspace{10pt}-\int^\tau_0\int_{\Omega}\langle Y^\epsilon_{t,x};\mathbf v\cdot\nabla_x\mathbf v\cdot\mathbf m+\rho^{(1)}\text{div}(q\mathbf v)\rangle dxdt
\end{align*}

By virtue of (4.3) and (6.2), we have
\begin{align*}
\text{div}_x(q\mathbf v)=\nabla_xq\cdot\mathbf v=\nabla_xq\cdot\nabla^{\perp}_x q=0
\end{align*}
and
\begin{align*}
-\int^\tau_0\int_{\Omega}\langle Y^\epsilon_{t,x};\partial_t\mathbf v\cdot\mathbf m+\partial_tq\rho^{(1)}\rangle dxdt
=-\int^\tau_0\int_{\Omega}\langle Y^\epsilon_{t,x};\omega\times\mathbf m\rangle\mathbf v\Delta_hq dxdt.
\end{align*}

Moreover, it is easy to check that
\begin{align*}
\mathbf v\cdot\nabla_x\mathbf v\cdot\mathbf m+(\omega\times\mathbf m)\cdot\mathbf v\Delta_hq=\mathbf m\cdot\nabla_h\frac{|\mathbf v|^2}{2}.
\end{align*}

So we deduce that
\begin{align*}
\frac{1}{2}\int^\tau_0\int_{\Omega}\langle Y^\epsilon_{t,x};\partial_t|\mathbf v|^2+\partial_t|q|^2\rangle dxdt
&=\frac{1}{2}\int^\tau_0\int_{\Omega}\langle Y^\epsilon_{t,x};\partial_t|\nabla_h^\perp q|^2+\partial_t|q|^2\rangle dxdt\\
&=\int^\tau_0\int_{\Omega}\langle Y^\epsilon_{t,x};\mathbf v\cdot\nabla_h(\Delta_h q)q\rangle dxdt\\
&=-\int^\tau_0\int_{\Omega}\langle Y^\epsilon_{t,x};\mathbf v\cdot\nabla_hq\Delta_hq\rangle dxdt=0,
\end{align*}
where the last equality due to $\mathbf v=\nabla^\perp q$.

Putting together Step 1 to Step 4, we conclude that
\begin{align*}
\mathcal{E}_\epsilon (\rho_\epsilon,\mathbf{m}_\epsilon|r_\epsilon, \mathbf U_\epsilon)+\mathcal{D}^{\epsilon}(\tau)\leq \omega(\epsilon)
+\int^\tau_0(1+\xi(t))[\mathcal{E}(\rho_\epsilon,\mathbf{m}_\epsilon|r_\epsilon, \mathbf U_\epsilon)+\mathcal{D}^\epsilon(t)]dt,
\end{align*}
where $r_\epsilon$, $\mathbf U_\epsilon$ are given by (6.1).
Letting $\epsilon\rightarrow0$ and applying the Gronwall's lemma, we complete the proof of Theorem 3.1.

\vskip 1.0cm

\vskip 0.5cm
\noindent {\bf Acknowledgements}

\vskip 0.1cm

The authors are grateful to the referee and the editor whose comments and suggestions greatly improved the presentation of this paper. The paper was written when Tong Tang was visiting the Institute of Mathematics of the Czech Academy of Sciences which {hospitality and support} is gladly acknowledged.


\end{document}